\documentclass[12pt,leqno]{article}
\usepackage{amsfonts}
\usepackage{amsmath, amsthm, amsfonts, amssymb, color}
\usepackage{mathrsfs}
\usepackage{color}
\theoremstyle{definition}

\newcommand{\scr}[1]{\mathscr #1}
\definecolor{wco}{rgb}{0.5,0.2,0.3}

\numberwithin{equation}{section} \theoremstyle{remark}
\newcommand{\ua}{\uparrow}

\title{{\bf Entropy Estimate for  Degenerate SDEs with Applications to Nonlinear Kinetic Fokker-Planck Equations}\footnote{ {\textbf{Funding: }} Feng-Yu Wang is supported in part by the National Key R\& D Program of China (No. 2022YFA1006000, 2020YFA0712900) and NNSFC (11921001). Panpan Ren is supported by NNSFC (12301180) and Research Centre for Nonlinear Analysis at Hong Kong PolyU.} }
\author{
{\bf Zhongmin Qian$^{b)}$\footnote{$^{b)}$ Mathematical Institute,  Oxford University,  Oxford, OX2 6GG, United Kingdom (zhongmin.qian@maths.ox.ac.uk). },   Panpan Ren$^{c)}$\footnote{$^{c)}$ Department of Mathematics, City University of  Hong Kong, Tat Chee Avenue, Hong Kong,  China (panparen@gmail.com).  },  Feng-Yu Wang$^{a)}$\footnote{  
$^{a)}$ Center for Applied Mathematics, Tianjin University, Tianjin 300072, China (wangfy@tju.edu.cn).}
}}

\begin{document}
\allowdisplaybreaks
\def\R{\mathbb R}  \def\ff{\frac} \def\ss{\sqrt} \def\B{\mathbf
B}
\def\N{\mathbb N} \def\kk{\kappa} \def\m{{\bf m}}
\def\ee{\varepsilon}\def\ddd{D^*}
\def\dd{\delta} \def\DD{\Delta} \def\vv{\varepsilon} \def\rr{\rho}
\def\<{\langle} \def\>{\rangle}
  \def\nn{\nabla} \def\pp{\partial} \def\E{\mathbb E}
\def\d{\text{\rm{d}}} \def\bb{\beta} \def\aa{\alpha} \def\D{\scr D}
  \def\si{\sigma} \def\ess{\text{\rm{ess}}}\def\s{{\bf s}}
\def\beg{\begin} \def\beq{\begin{equation}}  \def\F{\scr F}
\def\Ric{\mathcal Ric} \def\Hess{\text{\rm{Hess}}}
\def\e{\text{\rm{e}}} \def\ua{\underline a} \def\OO{\Omega}  \def\oo{\omega}
 \def\tt{\tilde}\def\[{\lfloor} \def\]{\rfloor}
\def\cut{\text{\rm{cut}}} \def\P{\mathbb P} \def\ifn{I_n(f^{\bigotimes n})}
\def\C{\scr C}      \def\aaa{\mathbf{r}}     \def\r{r}
\def\gap{\text{\rm{gap}}} \def\prr{\pi_{{\bf m},\varrho}}  \def\r{\mathbf r}
\def\Z{\mathbb Z} \def\vrr{\varrho} \def\ll{\lambda}
\def\L{\scr L}\def\Tt{\tt} \def\TT{\tt}\def\II{\mathbb I}
\def\i{{\rm in}}\def\Sect{{\rm Sect}}  \def\H{\mathbb H}
\def\M{\mathbb M}\def\Q{\mathbb Q} \def\texto{\text{o}} \def\LL{\Lambda}
\def\Rank{{\rm Rank}} \def\B{\scr B} \def\i{{\rm i}} \def\HR{\hat{\R}^d}
\def\to{\rightarrow}\def\l{\ell}\def\iint{\int}\def\gg{\gamma}
\def\EE{\scr E} \def\W{\mathbb W}
\def\A{\scr A} \def\Lip{{\rm Lip}}\def\S{\mathbb S}
\def\BB{\scr B}\def\Ent{{\rm Ent}} \def\i{{\rm i}}\def\itparallel{{\it\parallel}}
\def\g{{\mathbf g}}\def\Sect{{\mathcal Sec}}\def\T{\mathbb T}\def\BB{{\bf B}}
\def\f{\mathbf f} \def\g{\mathbf g}\def\BL{{\bf L}}  \def\BG{{\mathbb G}}
\def\Bd{{D^E}} \def\BdP{D^E_\phi} \def\Bdd{{\bf \dd}} \def\Bs{{\bf s}} \def\GA{\scr A}
\def\Bg{{\bf g}}  \def\Bdd{\psi} \def\supp{{\rm supp}}\def\div{{\rm div}}
\def\ddiv{{\rm div}}\def\osc{{\bf osc}}\def\1{{\bf 1}}\def\BD{\mathbb D}\def\GG{\Gamma}
\def\H{{\bf H}}

\maketitle

{\textbf{Abstract.}} The  relative entropy  for two different degenerate diffusion processes is estimated by using the Wasserstein distance of initial distributions and the difference between coefficients. As applications,  the entropy-cost inequality and exponential ergodicity in entropy are derived for distribution dependent stochastic  Hamiltonian systems associated with  nonlinear kinetic Fokker-Planck equations. \\

   \noindent
   

 \section{Introduction}

 To characterize  the stability of stochastic systems under perturbations, a natural way is to estimate the difference of  distributions for two different processes, see \cite{Qian} for a comparison theorem  on transition densities (i.e. heat kernels) of
  diffusions with different drifts.

 Recently, by using the entropy inequality  established  by  Bogachev,  R\"{o}ckner and  Shaposhnikov  \cite{BRS} for  diffusion processes,  and by developing a bi-coupling argument, the entropy and probability distances  have been estimated   in \cite{RW23, 23HRW} for different non-degenerate SDEs with distribution dependent noise.
 In this paper, we aim to establish entropy inequality for degenerate  diffusion processes. As applications, we establish a  log-Harnack inequality and study the exponential ergodicity in entropy for stochastic Hamiltonian systems with distribution dependent noise.

Let us start with a simple stochastic Hamiltonian system whose Hamiltonian function is given by
 $$H(x):= V_1(x^{(1)})+ V_2(x^{(2)})\ \ \text{for } x=(x^{(1)}, x^{(2)})\in\R^d\times\R^d,$$
 where $V_i \in C^2(\R^d)$ with $\|\nn^2 V_i\|_\infty<\infty, i=1,2$. Then $X_t=(X_t^{(1)},X_t^{(2)})$,
 the speed $X_t^{(1)}$ and the location $X_t^{(2)}$ of the stochastic particle, solves the following degenerate stochastic differential equation (SDE) on $\R^d\times\R^d:$
  \beq\label{1.1} \beg{cases} \d X_t^{(1)}= \nn V_2(X_t^{(2)})\d t,\\
\d X_t^{(2)}=  \ss 2\,  \d W_t -\big(\nabla V_1(X_t^{(1)})+\nn V_2(X_t^{(2)})\big)\d t,\end{cases}\end{equation}
where $W_t$ is the $d$-dimensional Brownian  motion on a filtered probability space $(\OO,\F,(\F_t)_{t\ge 0}, \P)$.
%$X_t^{(1)}$ and $X_t^{(2)}$ stand  for the location and speed at time $t$  for   a  random particle, respectively.
It is well known that the distribution density function of $X_t$ solves the associated  kinetic Fokker-Planck  equation.

When  for each $i=1,2$,  $ \mu^{(i)}(\d x^{(i)}):=\e^{-V_i(x^{(i)})}\d x^{(i)}$ is a probability measure on $\R^d$,   SDE \eqref{1.1} has a unique invariant probability measure
$$\bar\mu(\d x):=\mu^{(1)}(\d x^{(1)}) \mu^{(2)} (\d x^{(2)}),\ \ \text{for } x=(x^{(1)},x^{(2)}) \in \R^d\times\R^d.$$
According to   Villani \cite{V09}, suppose that $\mu^{(i)}$
 satisfies the Poincar\'e inequality
$$\mu^{(i)}(f^2)\le \mu^{(i)}(f)^2+ C \mu^{(i)}(|\nn f |^2),\ \ \forall f\in C_b^1(\R^d),i=1,2,$$
for some constant $C>0$, where and in the sequel $\mu(f):=\int f\d\mu$ for a measure $\mu$ and a function $f$ if the integral  exists.  Then the Markov semigroup $P_t$ associated with \eqref{1.1} converges exponentially to $\bar\mu$ in $H^1(\bar\mu)$, i.e. for some constants $c,\ll>0$,
$$
\bar\mu\big(|P_tf-\bar\mu(f)|^2+ |\nn P_tf|^2\big)\le c  \e^{-\ll t} \bar\mu\big(|f-\bar\mu(f)|^2+|\nn f|^2\big)
$$
for any $t\ge 0 $ and $f\in C_b^1(\R^d)$. This property, known as $``$hypocoercivity " due to Villani \cite{V09},  has been explored further by various authors in a series of papers for the exponential convergence of $P_t$ in $L^2(\mu)$,  such as  \cite{CHSG}  by  Camrud,  Herzog,   Stoltz and Gordina, as well as   \cite{GS} by  Grothaus and  Stilgenbauer,    based on an abstract analytic framework  built up  by Dolbeaut, Mouhot and Schmeiser \cite{DMS}, see also the recent work  \cite{GRW23} for the study of singular models.
In case the Poincar\'e inequality fails, slower convergence rates  are presented in   \cite{GW19, HW22} using the weak Poincar\'e inequality developed by R\"ockner and the third named author  \cite{RW01}.

On the other hand,
the study of the exponential ergodicity in the relative entropy arising from information theory, which is stronger than that in $L^2$ (see \cite{W17}), becomes an important topic. Recall that if  $\mu$ and $\nu$ are two probability measures, then the relative entropy of $\mu$ with respect to $\nu$ is defined by
$$\Ent(\mu|\nu):=\beg{cases} \mu\big(\log\ff{\d\mu}{\d\nu}\big),\ &\text{if} \ \mu\ \text{is\ absolutely\ continuous\ w.r.t.\ }\nu,\\
\infty,\ &\text{otherwise.}\end{cases}$$
By Young's inequality, see for instance \cite[Lemma 2.4]{09SPA}, for any positive measurable function $f$ such that $\nu(f)=1$, we have
$$\mu(\log f) = \nu\Big(\ff{\d\mu}{\d\nu} \log f\Big)\le \nu\Big(\ff{\d\mu}{\d\nu}\log\ff{\d\mu}{\d\nu}\Big)
+\log\nu(f)= \Ent(\mu|\nu),$$
and the equality holds for $f=\ff{\d\mu}{\d\nu}$. Thus,
\beq\label{*WY0} \Ent(\mu|\nu)= \sup_{f> 0, \nu(f)=1} \mu(\log f)= \sup_{f> 0, \nu(f)<\infty}\big[\mu(\log f)-\log \nu(f)\big],\end{equation}
since the right hand side is infinite if $\mu$ is not absolutely continuous with respect to $\nu$.

By establishing a log-Harnack inequality, the exponential ergodicity in entropy has been been
derived in  \cite{W17} for stochastic Hamiltonian systems for linear $\nn V_2$, and has been further extended in \cite{RW21, HW} to the case with distribution dependent drift.
However, the log-Harnack inequality and the exponential ergodicity in entropy are still unknown for  stochastic Hamiltonian systems with  nonlinear $\nn V_2$.

To formulate distribution dependent SDEs, we introduce the Wasserstein space $\scr P_2(\R^d)$ for probability measures on $\R^d$ having finite second moment. It is a Polish space   under the Wasserstein distance
$$\W_2(\mu,\nu):= \inf_{\pi\in\C(\mu,\nu)} \bigg(\int_{\R^{d}\times\R^{d}} |x-y|^2 \pi(\d x,\d y)\bigg)^{\ff 1 2},$$
where $\C(\mu,\nu)$ denotes the set of all couplings for $\mu$ and $\nu$.    Let $\L_\xi$ denote the distribution of the random variable $\xi$.

 To illustrate our general results,  we consider below the distribution dependent  stochastic Hamiltonian system for $X_t:=(X_t^{(1)},  X_t^{(2)})\in \R^{d_1}\times\R^{d_2}$:
 \beq\label{1.2} \beg{cases} \d X_t^{(1)}= \big\{ B X_t^{(2)} + b(X_t)\big\}\d t,\\
\d X_t^{(2)}=  \si(\L_{X_t})  \d W_t  +Z(X_t^{(2)},\L_{X_t}) \d t,\ \ t\ge 0,\end{cases}\end{equation}
where     $B$ is  a   $d_1\times d_2$-matrix such that $BB^*$ is invertible (i.e. Rank$(B)= d_1$), $b\in C_b^2(\R^{d_1+d_2})$ such that
$$\big\<(\nn^{(2)} b)B^*v,v\big\>\ge -\dd |B^*v|^2,\ \ \ v\in\R^{d_1} $$
 holds for some constant $\dd\in (0,1)$, where $\nn^{(2)}$ is the gradient in $x^{(2)}\in \R^{d_2},$
  and
$$
\si: \scr P_2(\R^{d_1+d_2} ) \to \R^{d_2\otimes d_2},\ \ Z: \R^{d_1+d_2}\times\scr P_2(\R^{d_1+d_2}) \to \R^{d_2}
$$
are Lipschitz continuous.  According to \cite[Theorem 2.1]{W18},  \eqref{1.2} is well-posed for distributions in $\scr P_2(\R^{d_1+d_2}),$ i.e. for any $\F_0$-measurable initial value $X_0$ with $\L_{X_0}\in \scr P_2(\R^{d_1+d_2})$, (respectively, any initial distribution $\mu\in \scr P_2(\R^{d_1+d_2})$), the SDE has a unique strong (respectively, weak) solution with $\L_{X_t}\in \scr P_2(\R^{d_1+d_2})$ continuous in $t\ge 0$.  Let $P_t^*\mu:=\L_{X_t}$ where $X_t$ is  the solution of (\ref{1.2}) with initial distribution $\mu\in \scr P_2$.
If $\nn Z(\cdot,\mu) $ is bounded and Lipschitz continuous uniformly in $\mu$,   then the following assertions are implied by Theorem 4.1.
\beg{enumerate} \item[$\bullet$] By \eqref{EN1*} for $k=0$,
 there exists a constant $c>0$ such that
$$
\Ent(P_t^*\mu|P_t^*\nu) \le \ff{c}{t^{3} }\W_2(\mu,\nu)^2,\ \ \ t\in (0,1];\ \mu,\nu\in \scr P_2(\R^{d_1+d_2}).$$
\item[$\bullet$] If $P_t^*$ is exponentially ergodic in $\W_2$, i.e. $P_t^*$ has a unique invariant probability measure $\bar\mu\in \scr P_2(\R^{d_1+d_2})$ and there exist two positive constants $c_1$ and $\ll$ such that
\beq\label{EX}
\W_2(P_t^*\mu, \bar\mu)^2\le c_1 \e^{-\ll t} \W_2(\mu,\bar\mu)^2
\end{equation}
holds for any $t\ge 0 $ and $\mu\in \scr P_2(\R^{d_1+d_2}) $,  then the exponential ergodicity in entropy holds:
$$
\Ent(P_t^*\mu|\bar\mu) \le cc_1\e^{-\ll(t-1)} \W_2(\mu,\bar\mu)^2
$$
holds for any $t\ge 0 $ and $\mu\in \scr P_2(\R^{d_1+d_2})$.
See Corollary \ref{CC} and Example 4.1 below for some concrete models satisfying \eqref{EX}.   \end{enumerate}

The remainder of the paper is organized as follows. We establish  an entropy inequality in Section 2 for some SDEs which applies also to the degenerate case,   then apply  the inequality  to stochastic Hamiltonian systems   and the distribution dependent model in Sections 3 and 4 respectively.

\section{Entropy estimate between diffusion processes}

Let $d,m\in\mathbb N, T\in (0,\infty)$, and $(W_t)_{t\in [0,T]}$ be   an $m$-dimensional Brownian motion on a filtered probability space $(\OO,\F,  (\F_t)_{t\in [0,T]},\P)$.  Consider the following SDEs on $\R^d$:
\beq\label{E0}
\d X_t^{\<i\>}= Z_i (t,X_t^{\<i\>})\d t+ \si_i(t,X_t^{\<i\>})\d W_t\ \  \text{ for }   t\in [0,T],
\end{equation}
where
$$
Z_i: [0,T]\times\R^d\to\R^d\ \text{ and } \si_i: [0,T]\times \R^d\to \R^{d\otimes m}
$$
are nice enough measurable maps   such that
 the SDE  is well-posed for $i=1,2$. Let
 $(P_{s,t}^{\<i\>})_{0\le s\le t\le T}$ be the corresponding Markov semigroups,  i.e.
 $$
 P_{s,t}^{\<i\>} f(x):= \E[f(X_{s,t}^{i,x})]\ \text{ for } f\in \B_b(\R^d) \text{ and } x\in \R^d,
 $$
 where $(X_{s,t}^{i,x})_{t\in [s,T]}$ solves \eqref{E0} for $t\in [s,T]$ with $X_{s,s}^{i,x}=x.$
 The corresponding generators are given by
$$
L_t^{\<i\>}:= {\rm tr}\big\{a_i(t,\cdot)\nn^2\big\} + Z_i(t,\cdot)\cdot\nn\ \ \text{ for }  t\in [0,T],
$$
where $a_i:=\ff 1 2 \si_i\si_i^*$ which may be degenerate.  If $v: [0,T]\mapsto\R^d$ is a path, then
$$
\|v\|_{a_2}(t) :=\sup_{x\in\R^d}\inf\big\{|w|: w\in\R^d, a_2(t,x)^{\ff 1 2}w=v(t)\big\}\ \
\text{ for } t\in [0,T],
$$
where the convention that $\inf\emptyset =\infty$ is applied.

Let $\scr P(\R^d)$ denote the space of all probability measures on $\R^d$.  For a given $\nu\in \scr P(\R^d)$,  $X_t^{i,\nu}$ denotes the solution to  \eqref{E0} with     $\L_{X_0^{i,\nu}}=\nu$, where and in the sequel, $\L_\xi$ stands for the law of a random variable $\xi$. Denote
$$
P_t^{i,\nu}=\L_{X_t^{i,\nu}} \ \ \text{ for } t\in [0,T],\ \nu\in \scr P (\R^d) \text{ and }  i=1,2.
$$
We shall make the following assumptions.

\beg{enumerate}
\item[$(A_1)$]
For any $0\le s\le t\le T,$ $P_{s,t}^{\<2\>} C_b^2(\R^d) \subset C_b^2(\R^d) $ so that  the   Kolmogorov backward equation holds for any $f\in C_b^2(\R^d)$:
$$
\pp_s P_{s,t}^{\<2\>}f= -L_s^{\<2\>} P_{s,t}^{\<2\>}f\ \ \text{ for } s\in [0,t] \text{ and } t\in (0,T].
$$
\item[$(A_2)$] For any $t\in (0,T]$, $(a_1-a_2)(t,\cdot)$ is differentiable on $\R^d$, and there exists
  a measurable function $H_{a_1-a_2}\cdot^{1,\nu}: (0,T]\mapsto (0,\infty)$ such that
\beg{align*}
&\big|\E \big[{\rm div }\{(a_1-a_2)(t,\cdot) \nn f\}(X_t^{1,\nu})\big]\big|\\
&\le H_{a_1-a_2}^{1,\nu}(t) \big(\E[|a_2(t,\cdot)^{\ff 1 2} \nn f|^2(X_t^{1,\nu})\big]\big)^{\ff 1 2}
\end{align*}
holds for any $ t\in (0,T]$ and $f\in C_b^2(\R^d)$.
\end{enumerate}

We remark that condition $(A_1)$ is satisfied   when the coefficients have bounded first and second order derivatives.  For the non-degenerate case,  it is satisfied for a class of H\"older continuous  $\si_2$ and $b_2$, see  for instance \cite{MPZ} and references within.
According to \cite{BRS},   condition $(A_2)$ is satisfied  if $a_2$ is invertible and $X_t^{1,\nu}$ has a distribution density $\rr_t^{1,\nu}$ such that $\log \rr_t^{1,\nu}$ is in a Sobolev space.
In this case, inequality \eqref{EI1} in the following theorem reduces to \cite[Theorem 1.1]{BRS}.
  In the next section,  we shall verify these conditions for some important examples of degenerate SDEs.

  We are now in a position to state and prove the main result.

\beg{thm}\label{T1} Assume that $(A_1)$ and $(A_2)$ are satisfied. Then
\beq\label{EI1} \Ent(P_t^{1,\nu}|P_t^{2,\nu})\le \ff 1 4\int_0^t \Big\{\|Z_1-Z_2-{\rm div}(a_1-a_2)\|_{a_2}(s) +H_{a_1-a_2}^{1,\nu}(s) \Big\}^2\d s\end{equation}
for any $t\in (0,T]$.
\end{thm}

\beg{proof}
Let $X_t^{i,\nu}$ solve \eqref{E0} with initial distribution $\nu$, and let $X_0^{1,\nu}=X_0^{2,\nu}$. Let $C_{b,+}^2(\R^d)$ denote the space of all functions  $f\in C_b^2(\R^d)$ such that $\inf f>0$. By \eqref{*WY0} and an approximation argument,
  we have
\beq\label{A0} \beg{split}&\Ent(P_t^{1,\nu}|P_t^{2,\nu})=\sup_{f\in C_{b,+}^2(\R^d)} I_t(f),\\
&I_t(f):= \E\log f(X_t^{1,\nu})- \log \E f(X_t^{2,\nu}).\end{split}\end{equation}
Noting that $(X_t^{2,x}: x\in \R^d)_{t\in [0,T]}$ is a (time inhomogenous) Markov process, for   any $f\in C_{b,+}^2(\R^{d}),$ we have
\beq\label{*WY} \E[f(X_t^{2,\nu})]= \int_{\R^d} (P_{0,t}^{\<2\>}f)\d\nu= \E[P_{0,t}^{\<2\>}f(X_0^{2,\nu})].\end{equation}
So, by Jensen's inequality, we obtain
\beq\label{A1}
\beg{split} I_t(f)&=\E\log f(X_t^{1,\nu})- \log \E (P_{0,t}^{\<2\>}f)(X_0^{2,\nu})\\
&\le  \E\log f(X_t^{1,\nu})- \E \log (P_{0,t}^{\<2\>}f)(X_0^{2,\nu})\\
&= \int_0^t \Big[\ff{\d }{\d s} \E \log (P_{s,t}^{\<2\>}f)(X_s^{1,\nu})\Big] \d s\end{split}
\end{equation}
for every $t\in (0,T]$.
By $(A_1)$
and using It\^o's formula for $X_s^{1,\nu},$   we derive that
\beg{align*} & \ff{\d }{\d s} \E \big(\log (P_{s,t}^{\<2\>}f)(X_s^{1,\nu})\big)= \E\bigg[\Big(L_s^{\<1\>}  \log (P_{s,t}^{\<2\>}f) -\ff{L_s^{\<2\>} P_{s,t}^{\<2\>}f}{P_{s,t}^{\<2\>}f}\Big)(X_s^{1,\nu})\bigg]\ \\
&= \E\Big[(L_s^{\<1\>}-L_s^{\<2\>}) \log (P_{s,t}^{\<2\>}f)(X_s^{1,\nu})-\big|\{a_2(s,\cdot)^{\ff 1 2}\nn \log P_{s,t}^{\<2\>}f\}\big|^2(X_s^{1,\nu})\Big]\\
&= \E\Big[{\rm div}\big\{(a_1-a_2)(s, \cdot)\nn \log P_{s,t}^{\<2\>} f\big\}(X_s^{1,\nu}) -\big|\{a_2(s,\cdot)^{\ff 1 2}\nn \log P_{s,t}^{\<2\>}f\big|^2(X_s^{1,\nu})\Big]\\
&\quad + \E\Big[\big\<\{Z_1-Z_2-{\rm div}(a_1-a_2)\}(s,\cdot),\nn\log P_{s,t}^{\<2\>}f\big\>(X_s^{1,\nu})\Big].\end{align*}
Combining this with $(A_2)$ gives  that
\beg{align*} & \ff{\d }{\d s} \E \big(\log (P_{s,t}^{\<2\>}f)(X_s^{1,\nu})\big)\\
&\le \Big[H_{a_1-a_2}^{1,\nu}(s) +\|Z_1-Z_2-{\rm div}(a_1-a_2)\|_{a_2}(s)\Big]\big(\E|a_2(s,\cdot)^{\ff 1 2}\nn
\log P_{s,t}^{\<2\>} f|^2(X_s^{1,\nu})\big)^{\ff 1 2}\\
&\quad -\E\Big[\big|a_2(s,\cdot)^{\ff 1 2}\nn \log P_{s,t}^{\<2\>}f\big|^2(X_s^{1,\nu})\Big]\\
&\le  \ff 1 4 \Big[H_{a_1-a_2}^{1,\nu}(s) +\|Z_1-Z_2-{\rm div}(a_1-a_2)\|_{a_2}(s)\Big]^2
\end{align*}
for every $s\in (0,t]$, which, together with \eqref{A0} and \eqref{A1}, implies  the desired estimate  \eqref{EI1}.
\end{proof}

As explained in \cite{RW23} that $|H_{a_1-a_2}^{1,\nu}(s) |^2$ is normally  singular for small $s$, such that   the upper bound in \eqref{EI1} becomes infinite. To derive a finite upper bound of the relative entropy,
we make use of the bi-coupling argument developed in \cite{RW23}, which leads to the following consequence where   different initial distributions are also allowed.

\beg{cor}\label{C2.1} Assume that $(A_1)$ and $(A_2)$ are satisfied, $H_{a_1-a_2}^{1,x}(s):= H_{a_1-a_2}^{1,\dd_x}(s)$ is measurable in $x\in\R^d$ such that
$$H_{a_1-a_2}^{1,\nu}=\int_{\R^d} H_{a_1-a_2}^{1,x}(s)\nu(\d x).$$
 Suppose that there exist a constant $p\in (1,\infty)$ and  a decreasing  function $\eta: (0,T]\mapsto (0,\infty)$ such that
\beq\label{HI}
|P_{s,t}^{\<2\>} f(x)|^p \le \big(P_{s,t}^{\<2\>} |f|^p(y) \big)\e^{\eta(t-s) |x-y|^2}
\end{equation}
for any $0\le s<t\le T$ and $f\in \B_b(\R^d)$.
Then there exists a constant $c>0$ such that
\beg{align*} \Ent(P_t^{1,\mu}|P_t^{2,\nu})\le &\,\inf_{\pi\in \C(\mu,\nu)} \int_{\R^d\times \R^d} \bigg(\ff p 4\int_{t_0}^t \Big\{\|b_1-b_2-{\rm div}(a_1-a_2)\|_{a_2}(s) +H_{a_1-a_2}^{1,x_1}(s) \Big\}^2\d s\\
&\,+ (p-1)\log \E\Big\{\exp\big[c \eta(t-t_0) \big|X_{t_0}^{1,x_1}-X_{t_0}^{2,x_2} \big|^2\big]\Big\}\bigg)\pi(\d x_1,\d x_2)\end{align*}
  for any $0<t_0<t\le T$ and $x,y\in\R^d$.
  \end{cor}

\beg{proof} For simplicity, denote $P_t^{i,x}= P_t^{i,\dd_x}$ where $i=1,2, x\in\R^d$, and $\dd_x$ is the Dirac measure at $x$.  Let $X_t(x_1)$ be the diffusion process starting from the initial value $x_1$  with the infinitesimal  generator given by
$$L_t:= 1_{[0,t_0]}(t) L_t^{\<1\>}+ 1_{(t_0,t]}(t) L_t^{\<2\>}.$$
Let $P_t^{\<t_0\>x_1}=\L_{X_t(x_1)}$.
By using \eqref{EI1} with $\nu=\dd_{x_1}$ and $P_t^{\<t_0\>x_1} $ in place of  $P_t^{2,x_1}$, and combining  with     \cite[(2.4) and (2.9)]{RW23}, we deduce that
\beq\label{PO} \beg{split} \Ent(P_t^{1,x_1}|P_t^{2,x_2})\le &\,\ff p 4\int_{t_0}^t \Big\{\|b_1-b_2-{\rm div}(a_1-a_2)\|_{a_2}(s) +H_{s}^{1,x_1}(a_1-a_2) \Big\}^2\d s\\
&\,+ (p-1)\log \E\bigg\{\exp\Big[c \eta(t-t_0) \big|X_{t_0}^{1,x_1}-X_{t_0}^{2,x_2} \big|^2\Big]\bigg\}.\end{split}\end{equation}
On the other hand, if $\pi\in \C(\mu,\nu)$, then  by using \eqref{A0},  \eqref{*WY}  and Jensen's inequality,  we obtain
\beg{align*}
&\Ent(P_t^{1,\mu}|P_t^{2,\nu}) =\sup_{f\in C_{b,+}^2(\R^d)} \big\{\E\log f(X_t^{1,\mu})- \log \E f(X_t^{2,\nu})\big\}\\
&= \sup_{f\in C_{b,+}^2(\R^d)} \bigg\{\int_{\R^d} P_t^{\<1\>} (\log f)(x_1) \mu(\d x_1)- \log \int_{\R^d} P_t^{\<2\>} f(x_2) \nu(\d x_2)  \bigg\}\\
&\le \sup_{f\in C_{b,+}^2(\R^d)} \bigg\{\int_{\R^d} P_t^{\<1\>} (\log f)(x_1) \mu(\d x_1)-  \int_{\R^d} \log P_t^{\<2\>} f(x_2) \nu(\d x_2)  \bigg\}\\
&= \sup_{f\in C_{b,+}^2(\R^d)}\int_{\R^d\times\R^d}  \big\{  P_t^{\<1\>} (\log f)(x_1)  -    \log P_t^{\<2\>} f(x_2)   \big\}\pi(\d x_1,\d x_2)\\
&\le \int_{\R^d\times\R^d}  \sup_{f\in C_{b,+}^2(\R^d)} \big\{  P_t^{\<1\>} (\log f)(x_1)  -    \log P_t^{\<2\>} f(x_2)    \big\} \pi(\d x_1,\d x_2)\\
&= \int_{\R^d\times\R^d}  \Ent( P_t^{1,x_1} | P_t^{2,x_2})     \pi(\d x_1,\d x_2),
\end{align*}
which, together with \eqref{PO}, yields the desired estimate.
\end{proof}

\section{Stochastic Hamilton system}

\subsection{A general result  }
Let $d_1,d_2\in\mathbb N$. For any  initial distribution $\nu\in \scr P(\R^{d_1+d_2})$, consider the following degenerate  SDEs for $X_t^{i,\nu}= (X_t^{i(1),\nu}, X_t^{i(2),\nu})\in \R^{d_1}\times \R^{d_2}$ ($i=1,2$):
\beq\label{EI}
\beg{cases}
\d X_t^{i(1),\nu}= \tt b (t,X_t^{i,\nu})\d t,\\
\d X_t^{i(2),\nu}=  Z_i(t,X_t^{i,\nu}) \d t
+ \si_i(t,X_t^{i,\nu})\d W_t,\ \ \ \L_{X_0^{i,\nu}}=\nu, \
\text{ for }t\in [0,T],
\end{cases}
\end{equation}
where $W_t$ is a $d_2$-dimensional Brownian motion on a filtered probability space $(\OO, \F, (\F_t)_{t\in [0,T]},\P)$, and
$$\tt b: [0,T]\times \R^{d_1+d_2}\to \R^{d_1},\ \ Z_i: [0,T]\times \R^{d_1+d_2}\to \R^{d_2},\ \ \si_i: [0,T]\times \R^{d_1+d_2}\to \R^{d_2\otimes d_2}$$ are measurable.

If $\nu=\dd_x$ where $x\in \R^{d_1+d_2}$,  then the solution is simply denoted by $X_t^{i,x}=(X_t^{i(1),x},X_t^{i(2),x})$.
Let $\nn^{(i)}$ be the gradient in $x^{(i)}\in \R^{d_i}$ for $i=1,2$.

Let us introduce the following technical conditions.

\beg{enumerate} \item[$(B_1)$] The coefficients $\si_i(t,x), Z_i(t,x)$ (for $i=1,2$)  and $\tt b(t,x)$ are locally bounded in $(t,x)\in [0,T]\times \R^{d_1+d_2}$ and   twice differentiable  in the space variable $x$. The matrix valued function $a_2:=\ff 1 2 \si_2\si_2^*$ is invertible.  There exists a constant $K>0$  such that
\beg{align*} \|\nn^j Z_i(t,x)\|+\|\nn^j \tt b(t,x)\|+\|\nn^j \si_i(t,x)\|\le K\end{align*}
for  $(t,x)\in [0,T]\times \R^{d_1+d_2}$ and $j=1,2$.
 \item[$(B_2)$] There exists a   function  $\xi^{\nu}\in C((0,T];(0,\infty))$ such that
$$ \big|\E [(\nn_v^{(2)} f)(X_t^{1,\nu})]\big| \le \xi_t^{\nu}\big(\E[ f(X_t^{1,\nu})^2]\big)^{\ff 1 2}
$$
for $t\in (0,T]$, $v^{(2)}\in \R^{d_2}$ with $|v^{(2)}|=1 $ and  $f\in C_b^1(\R^{d_1+d_2})$.
\end{enumerate}

It is well known  that condition $(B_1)$ implies the well-posededness of \eqref{EI} and that condition $(A_1)$ is satisfied.
Let  $P_t^{i,\nu}$ be the distribution of $X_t^{i,\nu}$.

To state our next result we recall that for a vector valued function $g$ on $[0,T]\times\R^{d_1+d_2}$
$$\|g\|_{t,\infty}:=\sup_{z\in \R^{d_1+d_2}} |g(t,z)|
$$
for $t\in [0,T]$.

\beg{thm}\label{T2}
Assume that conditions $(B_1)$ and $(B_2)$ are satisfied. Let $(e_j)_{1\le j\le d_2}$ be the canonical basis on $\R^{d_2}$.

1) The following equality holds:
\beg{align*}
&\Ent(P_t^{1,\nu}|P_t^{2,\nu})  \le \ff 1 4 \int_0^t \Big[  \big\|a_2^{-\ff 1 2} \big\{Z_1-Z_2-{\rm div}(a_1-a_2)\big\} \big\|_{s,\infty}
+ \xi_s^{\nu}\sum_{j=1}^{d_2}   \big\| a_2^{-\ff 1 2} (a_1-a_2)e_j\big\|_{s,\infty} \Big]^2\d s.
\end{align*}

2) Suppose $\eqref{HI}$ holds, then there exists a constant $c>0$ such that
\beg{align*}
\Ent(P_t^{1,\mu}|P_t^{2,\nu})\le &\,\inf_{\pi\in \C(\mu,\nu)} \int_{\R^{d_1+d_2}\times \R^{d_1+d_2} }\bigg(  p  I_{t_0,t}^{x_2}
 + (p-1)\log \E\Big[\e^{c \eta(t-t_0) |X_{t_0}^{1,x_1}-X_{t_0}^{2,x_2} |^2} \Big] \bigg)\pi(\d x_1,\d x_2)
 \end{align*}
 for any $0<t_0<t\le T$ and  $\mu,\nu\in \scr P(\R^{d_1+d_2})$,
where
$$
I_{t_0,t}^x:= \ff 1 4 \int_{t_0}^t \Big[  \big\|a_2^{-\ff 1 2} \big\{Z_1-Z_2-{\rm div}(a_1-a_2)\big\} \big\|_{s,\infty}
+ \xi_s^x\sum_{j=1}^{d_2}   \big\| a_2^{-\ff 1 2} (a_1-a_2)e_j\big\|_{s,\infty} \Big]^2\d s
$$
and $\xi^x_s:=\xi_s^{\dd_x}$ for every   $x\in \R^{d_1+d_2}$ and $s\in [t_0,t]$.
\end{thm}

\beg{proof}  As explained in the proof of Corollary \ref{C2.1}, we only need to prove the first estimate.
Since $(B_2)$ is satisfied, we have
 \beg{align*} &\big|\E \big[{\rm div }\{{\rm diag}\{{\bf 0}_{d_1\times d_1},(a_1-a_2)(t,\cdot)\} \nn f\}(X_t^{1,\nu})\big]\big|\\
&=\Big|\sum_{j=1}^{d_2}  \E \big[\pp_{y_j} \{(a_1-a_2)(t,\cdot) \nn^{(2)} f\}_j\big](X_t^{1,\nu})\Big|\\
&\le \xi^{\nu}_t\sum_{j=1}^{d_2} \big(\E \{(a_1-a_2)(t,\cdot)\nn^{(2)} f\}_j(X_t^{1,\nu})^2 \big)^{\ff 1 2}\\
& = \xi^{\nu}_t\sum_{j=1}^{d_2} \big(\E \<a_2(t,\cdot)^{-\ff 1 2}(a_2-a_2)(t,\cdot) e_j, a_2(t,\cdot)^{\ff 1 2} \nn^{(2)}f\>_{\R^{d_2}}(X_t^{1,\nu})^2 \big)^{\ff 1 2}\\
&\le \xi^{\nu}_t \sum_{j=1}^d   \| |a_2^{-\ff 1 2} (a_1-a_2)e_j|\|_{t,\infty} \big(\E\big|a_2(t,\cdot)^{\ff 1 2} \nn^{(2)} f\big|^2(X_t^{1,\nu})\big)^{\ff 1 2}.\end{align*}
Thus  $(A_2)$ is satisfied  with
$$H_t^\nu(a_1-a_2):= \xi^{\nu}_t  \sum_{j=1}^d   \| |a_2^{-\ff 1 2} (a_1-a_2)e_j|\|_{t,\infty}.$$
Since $(B_1)$ implies $(A_1)$, the desired estimate follows immediately from Theorem \ref{T1}.
\end{proof}

\subsection{A class of   models}
We next discuss a class of degenerate stochastic models for which
condition $(B_2)$ is satisfied and the dimension-free Harnack inequality   \eqref{HI} holds.

Consider the following  SDE for $X_t^{i,\nu}=(X_t^{i(1),\nu}, X_t^{i(2),\nu})\in \R^{d_1+d_2}$:
 \beq\label{1.2'}
 \beg{cases}
 \d X_t^{i(1),\nu}= \big\{A X_t^{i(1),\nu}+ B X_t^{i(2),\nu}+b(X_t^{i,\nu})\big\}\d t,\\
\d X_t^{i(2),\nu}=  \si_i(t)  \d W_t  +Z_i(t,X_t^{i,\nu}) \d t,\ \ \L_{X_0^{i,\nu}}=\nu \ \text{ for } i=1,2,
\end{cases}
\end{equation}
where   $A$, $B$,  $b$, $\si_i$ and $Z_i$ satisfy the following assumption.

\beg{enumerate} \item[$(B_3)$]
1) $A$ is a $d_1\times d_1$ matrix and $B$ is a $d_1\times d_2$ matrix, such that Kalman's condition
\beq\label{R}
{\rm Rank}\left[A^iB: 0\le i\le k\right]=d_1
\end{equation}
holds for some $0\le k\le d_1-1$.

2) $b\in C_b^1(\R^{d_1+d_2})$ with Lipschitz continuous $\nn b$, and there exists a constant $\dd\in (0,1)$ such that
\beq\label{R'} \big\<(\nn^{(2)}b(x))B^*v, v\big\>\ge -\dd |B^*v|^2,\ \ v\in\R^{d_1}, x\in \R^{d_1+d_2}.\end{equation}

3) $\si_1(t)$ and $\si_2(t)$ are  bounded, and $a_2(t):= \ff 1 2 \si_2(t)\si_2(t)^*$ is invertible with bounded inverse.

4) $Z_i(t,x)$ (for $i=1,2$) are locally bounded in $[0,T]\times\R^{d_1+d_2}$ and differentiable in $x$,  such that
$$
\sup_{t\in [0,T]} \bigg\{\|\nn Z_i(t,\cdot)\|+ \ff{\|\nn Z_i(t,x) -\nn Z_i(t,y)\|}{|x-y|}\bigg\}\le K
$$
holds for some constant $K>0$.
\end{enumerate}
We introduce $\xi_t$ in two different cases:
\beq\label{XI'}    \xi_t:=\beg{cases}
 t^{-2k-\ff 1 2}, &\text{if}\  Z_1(t,x)= Z_1(t,x^{(2)})\ \text{\ is\ independent\ of}\ x^{(1)}, \\
  t^{-2k-\ff 3 2}, &\text{otherwise}. \end{cases}  \end{equation}

\beg{cor}\label{C3.1}
Assume that $(B_3)$ is satisfied for either $k=0$ or $k\ge 1$ but $b(x)=b(x^{(2)})$ only depends on $x^{(2)}$.  Let $P_t^{i,\nu}$ be the distribution of $X_t^{i,\nu}$ solving $\eqref{1.2'}$. Then there exist    constants $c>0$ and $\vv\in (0,\ff 1 2]$ such that for any $ t\in (0,T]$ and $\mu,\nu\in \scr P_2(\R^{d_1+d_2}),$
\beg{align*}
\Ent(P_t^{1,\nu}|P_t^{2,\mu})\le   \ff c {t^{4k+3}} \bigg(\W_2(\mu,\nu)^2 + \int_0^t \|Z_1-Z_2\|_{s,\infty}^2\bigg)
 + c\int_{\vv(1\land t)^{4k+3}}^t  \xi_s^2  \|a_1(s)-a_2(s)\|^2\big)\d s.
\end{align*}
\end{cor}

\beg{proof}
Without loss of generality, we may and do assume that $  \si_i=  \ss{2 a_i}.$ Moreover, by a standard approximation argument, under $(B_3)$ we may find a sequence
$\{Z_i^{(n)}\}_{n\ge 1}$ for each $i=1,2$, such that
\beg{align*}&\sup_{n\ge 1,k=1,2, t\in [0,T]} \|\nn^k Z_i^{(n)}(t,\cdot)\|\le K,\\
& \lim_{n\to\infty}\sup_{t\in [0,T]} \big\{\|(Z_i-Z_i^{(n)})(t,\cdot)\|_\infty+\|\nn(Z_i-  Z_i^{(n)})(t,\cdot)\|_\infty=0.\end{align*}
Moreover, let $\{b^{(n)}\}_{n\ge 1}$ be a bounded sequence in $ C_b^2(\R^{d_1+d_2})$   such that $\|b^{(n)}-b\|_{C_b^1(\R^{d_1+b_2})}\to 0$ as $n\to\infty.$
Let  $P_t^{i,\nu,;n}$ be defined as $P_t^{i,\nu}$ for $(b^{(n)},Z_i^{(n)})$ replacing $(b,Z_i).$  It is well known that $P_t^{i,\nu;n}\to P_t^{i,\nu}$ weakly as $n\to\infty$, so that \eqref{A0}
implies that
 $$\Ent(P_t^{1,\nu}|P_t^{2,\mu})\le \liminf_{n\to\infty} \Ent(P_t^{1,\nu;n}|P_t^{2,\mu;n}).$$
  Therefore, we may and do assume that $\|\nn^k b\|+\|\nn^k Z_i(t,\cdot)\|_\infty\le K$ holds for some constant $K>0$ and $i,k=1,2$, so that Theorem \ref{T2} applies.

(a) By $(B_3)$,   $\si_1\ge 0, \si_2\ge \ll I_{d_2}$ for some constant $\ll>0$, where $I_{d_2}$ is the $d_2\times d_2$ identity  matrix. So,  according to the proof of \cite[Lemma 3.3]{PW06},
\beq\label{AS}   \|\si_1-\si_2\|=\bigg\|2\int_0^\infty\e^{-r \si_1}(a_1-a_2) \e^{-r \si_2}\d r  \bigg\|\le \ff 2 \ll \|a_1-a_2\|.\end{equation}
By Lemma \ref{TN1} below,  there exists a constant $c_1>0$ such that for any $\nu$, condition $(B_2)$ holds with
\beq\label{XI}  \xi^\nu_t= c_1 \xi_t:=\beg{cases} c_1
t^{-2k-\ff 1 2}, &\text{if}\ Z_1(t,x)= Z_1(t,x^{(2)}),\\
c_1 t^{-2k-\ff 3 2}, &\text{in\ general}. \end{cases}  \end{equation}
Moreover, by Lemma \ref{LNN} below,
\eqref{HI} holds for the following $\eta(s), s\in (0,T):$
\beq\label{ETA} \eta(s)= c(p)
s^{-4k-3}, \  \ s\in (0,T].  \end{equation}   Combining these with Theorem \ref{T2}, and noting that $a_2^{-1}$ is bounded and ${\rm div}(a_1-a_2)=0$, we can find a constant $c_2>0$ such that
for any $0<t_0<t\le T,$
\beq\label{YY1} \beg{split}& \Ent(P_t^{1,\mu}|P_t^{2,\nu})\le c_2 \int_{t_0}^t  \Big(\big\|Z_1-Z_2\big\|_{s,\infty} ^2 +  |\xi_s|^2  \big\|a_1(s)-a_2(s)\big\|^2\Big)\d s\\
& + c_2 \inf_{\pi\in \C(\mu,\nu)}\int_{\R^{d_1+d_2}\times\R^{d_1+d_2} } \log \E\Big[\e^{c_2 (t-t_0)^{-4k-3}  |X_{t_0}^{1,x_1}-X_{t_0}^{2,x_2} |^2} \Big]  \pi(\d x_1,\d x_2).\end{split}\end{equation}
It remains to estimate the exponential expectation in the last term.

(b) By $(B_3)$ and \eqref{AS}, there exists a constant $c_3\ge 1$ such that
$$
\d |X_s^{1,x_1}-X_s^{2,x_2}|^2\le c_3\big(|X_s^{1,x_1}-X_s^{2,x_2}|^2+ \|Z_1-Z_2\|_{s,\infty}^2 +\|a_1(s)-a_2(s)\|^2\big)\d s + \d M_s,
$$
where
$$\d M_s:= 2\<X_s^{1,x_1}-X_s^{2,x_2}, \{\si_1(s)-\si_2(s)\}\d W_s\>
$$
and therefore the following differential inequality holds:
\beq\label{PP1} \d \<M\>_s\le c_3 |X_s^{1,x_1}-X_s^{2,x_2}|^2\d s.\end{equation}
It follows that
\beq\label{PP2}
\beg{split} & |X_s^{1,x_1}-X_s^{2,x_2}|^2\le \e^{c_3 s} |x_1-x_2|^2\\
&\quad  + \int_0^s \e^{c_3(s-r)} \big(\|Z_1-Z_2\|_{r,\infty}^2 +\|a_1(r)-a_2(r)\|^2\big)\d r+\int_0^s  \e^{c_3(s-r)} \d M_r.\end{split}
\end{equation}
Let
$$
\tau_n:=\inf\big\{s\in [0,T]: |X_s^{1,x_1}-X_s^{2,x_2}|\ge n\big\},\ \text{ for } n= 1,2,\cdots
$$
with the convention that $\inf \emptyset:=T$. Then $\tau_n\to T$ as $n\to\infty$. Let
$$\ll:= c_3(t-t_0)^{-4k-3},\ \ c_4   :=\e^{c_3 T}.$$
By \eqref{PP2} and the fact  that 
$$\E[\e^{\ll\hat{N}_{t\land \tau_n}}]\le (\E\e^{2\ll^2\<\hat{N}\>_{t\land\tau_n}})^{\ff 1 2}\le (\E\e^{2\ll^2c_4^2\<M\>_{t\land\tau_n}})^{\ff 1 2},\ \ \ll\ge 0$$ holds for the continuous martingale
$$
\hat{N}_{t}:=\int_0^{t} \e^{c_3(s-r)}\d M_r,  \ \  \ t\ge 0,  
$$
we deduce that
\beq\label{*D1} \beg{split}&\E\big[\e^{\ll |X_{s\land\tau_n}^{1,x_1}-X_{s\land\tau_n}^{2,x_2}|^2}\big] \\
&\le \e^{c_4\ll |x_1-x_2|^2 + c_4\ll \int_0^s  \big(\|Z_1-Z_2\|_{r,\infty}^2 +\|a_1(r)-a_2(r)\|^2\big)\d r}
\big(\E\big[\e^{2\ll^2c_4^2\<M\>_{s\land\tau_n}}\big]\big)^{\ff 1 2}. \end{split}\end{equation}
While by \eqref{PP1} and Jensen's inequality,
\beq\label{*D2} \beg{split} &\E\big[\e^{2\ll^2c_4^2\<M\>_{s\land\tau_n}}\big]\le \E\Big[\e^{2\ll^2 c_4^2 c_3^2 \int_0^s |X_{r\land\tau_n}^{1,x_1}-X_{r\land\tau_n}^{2,x_2}|^2\d r}\Big]\\
&\le \ff 1 s\int_0^s \E\big[\e^{2\ll^2 c_4^2 c_3^2 s   |X_{r\land\tau_n}^{1,x_1}-X_{r\land\tau_n}^{2,x_2}|^2}\big]\d r \\
&\le \sup_{r\in [0,t_0]} \E\big[\e^{2\ll^2 c_4^2 c_3^2 t_0  |X_{r\land\tau_n}^{1,x_1}-X_{r\land\tau_n}^{2,x_2}|^2}\big] \end{split} \end{equation}
for $ s\in [0,t_0]$. Choosing
\beq\label{*D3} t_0 =\ff 1 {2 c_4^2 c_3^3}\Big(\ff{1\land t}2\Big)^{4k+3}=:\vv (1\land t)^{4k+3}\end{equation}
such that
$$2\ll c_4^2c_3^2 t_0 = 2 c_4^2c_3^3(t-t_0)^{-4k-3}t_0 \le 1,$$
we therefore conclude from \eqref{*D1} and \eqref{*D2} that
\beg{align*}  &\sup_{s\in [0,t_0]} \E\big[\e^{\ll |X_{s\land\tau_n}^{1,x_1}-X_{s\land\tau_n}^{2,x_2}|^2}\big] \\
 &\le \e^{c_4\ll |x_1-x_2|^2 + c_4\ll \int_0^{t_0}  \big(\|Z_1-Z_2\|_{r,\infty}^2 +\|a_1(r)-a_2(r)\|^2\big)\d r}
   \Big(\sup_{s\in [0,t_0]} \E\big[\e^{\ll |X_{s\land\tau_n}^{1,x_1}-X_{s\land\tau_n}^{2,x_2}|^2}\big]\Big)^{\ff 1 2}.\end{align*}
This together with the definition of $\ll$ and Fatou's lemma yields
\beg{align*} &   \E\big[\e^{c_3 (t-t_0)^{-4k-3}  |X_{t_0}^{1,x_1}-X_{t_0}^{2,x_2}|^2}\big]\le  \liminf_{n\to\infty}  \E\big[\e^{\ll |X_{t_0\land\tau_n}^{1,x_1}-X_{t_0\land\tau_n}^{2,x_2}|^2}\big] \\
&\le \e^{2c_4\ll |x_1-x_2|^2 + 2c_4\ll \int_0^{t_0}  \big(\|Z_1-Z_2\|_{r,\infty}^2 +\|a_1(r)-a_2(r)\|^2\big)\d r}.\end{align*}
Combining \eqref{YY1} with  \eqref{*D3}, we can therefore find a constant $c_5>0$ such that
\beg{align*} \Ent(P_t^{1,\mu}|P_t^{2,\nu})\le & c_2 \int_{\vv(1\land t)^{4k+3}}^t \Big( \big\|Z_1-Z_2\big\|_{s,\infty} ^2 +  |\xi_s|^2  \big\|a_1(s)-a_2(s)\big\|^2\Big)\d s\\
& + \ff{c_5}{t^{4k+3} } \bigg(\W_2(\mu,\nu)^2+  \int_0^{\vv(t\land 1)^{4k+3}}    \big(\|Z_1-Z_2\|_{r,\infty}^2 +\|a_1(r)-a_2(r)\|^2\big)\d r\bigg).\end{align*}
The desired estimate   now follow from \eqref{XI} immediately.
\end{proof}

\subsection{ Verify  conditions $(B_2)$ and \eqref{HI}}

Let us consider   $X_t=(X_t^{(1)}, X_t^{(2)})$ taking values in $\R^{d_1}\times \R^{d_2}$, which solves the SDE:
\beq\label{DS}
\beg{cases}
\d X_t^{(1)}= \big\{AX_t^{(1)}   + B X_t^{(2)} + b(X_t)\big\} \d t,\\
\d X_t^{(2)}=  Z(t,X_t) \d t
+ \si(t)\d W_t\ \ \text{for } t\in [0,T].
\end{cases}
\end{equation}
We have the following result which ensures condition $(B_2)$.

\beg{lem}\label{TN1} Let  $A,B,b$ and $(Z_i,\si_i):=  (Z,\si)$  satisfy conditions in $(B_3)$, but $b$ is not necessarily bounded.   Let $\xi_t$ be in $\eqref{XI'}.$  Then for any $p>1$ there exists a constant $c(p)>0$ such that  for any solution $X_t$ of \eqref{DS},
\beq\label{ES0}
\sup_{v\in \R^{d_1+d_2}, |v|=1}\big|\E\big[(\nn_{v}  f)(X_t)\big]\big|\le c(p) t^{-2 k-\ff 3 2}  \big(\E |f(X_t)|^p\big)^{\ff 1 p},\ \  t\in (0,T],  f\in C_b^1(\R^{d_1+d_2}).
\end{equation}    If $Z (t,x)= Z^ (t,x^{(2)})$ does not depend on $x^{(1)}$, then
\beq\label{ES}
\sup_{v\in \R^{d_2}, |v|=1}\big|\E\big[(\nn_{v}^{(2)} f)(X_t)\big]\big|\le c(p) t^{-2k-\ff 1 2} \big(\E |f(X_t)|^p\big)^{\ff 1 p},\ \   t\in (0,T],  f\in C_b^1(\R^{d_1+d_2}).
\end{equation}

  \end{lem}

\beg{proof}  We will follow the line of \cite[Remark 2.1]{WZ13} to establish the integration by parts formula
$$\E\big[(\nn_v f)(X_t)\big]=\E\big[f(X_t) M_t\big]$$
for some random variable $M_t\in L^{\ff p{p-1}}(\P).$ To this end, we first estimate $D_h X_t$ and $D_h (\nn X_t)^{-1}$, where $D_h$ is the Malliavin derivative along   an adapted process $(h_s)_{s\in [0,t]}$ on $\R^d$ with $$\E\int_0^t |h_s'|^2\d s<\infty.$$

(a)  For any $s\in [0,T)$, let $\{K(t,s)\}_{t\in [s,T]}$ solve the following random ordinary differential equation on $\R^{d_1\otimes d_1}$:
$$
\pp_t K_{t,s}= \big\{A X_t^{(1)}+ \nn^{(1)} b(t, X_t) \big\} K_{t,s},
\ \  K_{s,s}=I_{d_1}\ \text{ for } t\in [s,T].
$$
  Since $\nn b$ is bounded,    $K_{t,s}$ is bounded and invertible satisfying
\beq\label{Q0} \|K_{t,s}\|\lor\|K_{t,s}^{-1}\|\le \e^{K(t-s)}\ \ \text{ for } 0\le s\le t\le T \end{equation}
for some constant $K>0$.

Let
$$
Q_{t,s}:= \int_0^s\ff{r(t-r)}{t^2} K_{t,r}BB^*K_{t,r}^*\d r
\ \text{ for } 0\le s\le t\le T.
$$
By \cite[Theorem 4.2(1)]{WZ13} for $(t,s)$ replacing $(T,t)$,  when $k\ge 1$ and $b(x)= b(x^{(2)})$,  conditions \eqref{R} and \eqref{R'}  imply  that
\beq\label{Q}
Q_{t,s}\ge \ff{ c_0} t  s^{2(k+1)}I_{d_1}=: \xi_{t,s}I_{d_1}\ \ \text{ for } 0<s\le t\le T
\end{equation}
holds for some constant $c_0>0$. It is easy to see that this estimate also holds for $k=0$ and bounded  $\nn b(x)$ since in this case $BB^*$ is invertible.

Let $X_t(x)=(X_t^j(x))_{1\le j\le d_1+d_2}$ be the solution to \eqref{DS} with $X_0(x)=x.$ Since $\nn b$ and $\nn Z$ are bounded,   we see that
$$\nn X_t(x):=(\pp_{x_i}X_t^j(x))_{1\le i,j\le d_1+d_2}$$ exists and is invertible, and   the inverse $(\nn X_t(x))^{-1}=\big((\nn X_t(x))^{-1}_{ki}\big)_{1\le k,i\le d_1+d_2}$ satisfies
\beq\label{Q2}
\big\|\{\nn X_t(x)\}^{-1}\big\|\le c_1\ \ \text{ for } t\in [0,T]
\end{equation}
 for some constant $c_1>0$.

(b) Since $\nn b$ and $\nn Z$ are bounded,  $(D_h X_s)_{s\in [0,t]}$  is the unique solution of the random ODE
$$
\beg{cases}
\pp_s \big\{D_h X_s^{(1)}\big\}= AD_{h}X_s^{(1)} + B D_{h}X_s^{(2)}+\nn_{D_h X_s} b( X_s),\\
\pp_s \big\{D_h X_s^{(2)}\big\} = \nn_{D_h X_s} Z  (s, X_s)+\si(s) h_s',\ \ \ D_h X_0=0 \ \ \text{ for } s\in [0,t],
\end{cases}
$$
and there exists a constant $c_2>0$ such that
\beq\label{Q3}
|D_h X_s|\le c_2 \int_0^s |h_r'|\d r \ \ \text{ for } s\in [0,t].
\end{equation}
Similarly, since $\nn^2 b$ and $\nn^2 Z$ are also bounded,  for any $v\in \R^{d_1+d_2},$ $(D_h \nn_v X_s)_{s\in [0,t]}$ solve the equations
$$\beg{cases}
\pp_s \big\{D_h \nn_v X_s^{(1)}\big\}=  A D_{h}\nn_v X_s^{(1)}+ B D_h \nn_v X_s^{(2)} + \nn_{D_h\nn_v X_s} b(X_s)\\
  \qquad\qquad\qquad \qquad + \Big\{\nn^2 b(X_s)\Big\}\Big(D_h X_s, \nn_v X_s\Big) \\
\pp_s \big\{D_h \nn_v X_s^{(2)}\big\}= \nn_{D_h \nn_v X_s} Z (s, X_s) + \big\{\nn^2 Z (s,X_s)\big\}\big(D_h X_s, \nn_v X_s\big)
\end{cases} $$ for $ D_h \nn_v X_0=0$ and $s\in [0,t].$ Moreover,  there exists a constant $c_3>0$ such that
\beq\label{Q4}
\sup_{v\in \R^{d_1+d_2}, |v|\le 1} \big\| D_h \nn_v X_t\big\|\le c_3 \int_0^t \d s \int_0^s |h_r'|\d r\le c_3  t \int_0^t |h_s'|\d s.
\end{equation}

(c) For any fixed $t\in (0,T]$, we may construct $h$ by means of \cite[(1.8) and (1.11)]{WZ13} for $t$ replacing $T$ with the specific choice $\phi(s):=\ff{s(t-s)}t$ satisfying $\phi(0)=\phi(t)=0$ as required therein.

For any  $v=(v^{(1)},v^{(2)})\in \R^{d_1}\times\R^{d_2},$  let
\beg{align*} &\aa_{t,s}(v):=\ff{t-s}t v^{(2)} -\ff{s(t-s)}{t^2} B^*K_{t,s}^*Q_{t,t}^{-1}\int_0^t \ff{t-r}t K_{t,r} Bv^{(2)} \d r \\
&\qquad \qquad -\ff{s(t-s)B^*K_{t,s}^*}{\xi_{t,s}^2\d s}\int_0^t \xi_{t,s}^2Q_{t,s}^{-1} K_{t,s} v^{(1)} \d s,\\
&g_{t,s}(v):= K_{s,0}v^{(1)} +\int_0^s K_{s,r}B \aa_{t,s}(v)\d s,\\
&h_{t,s}(v):= \int_0^s \si(r)^{-1} \big\{\nn_{(g_{t,r}(v), \aa_{t,r}(v))} b (r, X_r)-\pp_r \aa_{t,r} \big\}\d r\ \ \text{ for } s\in [0,t].\end{align*}
 Let $\{e_i\}_{1\le i\le d_1+d_2}$ be the canonical ONB on $\R^{d_1+d_2}.$  According to \cite[Remark 2.1]{WZ13}, we have
\beq\label{Q5}\beg{split} & \E\big[(\nn_{e_i} f)(X_t\big]= \E\big[f(X_t) M_t(e_i)\big],\\
& M_t(e_i):=\sum_{j=1}^{d_1+d_2} \Big\{\dd(h_{t,\cdot}(e_j))(\nn X_t)^{-1}_{ji}-D_{h_{t,\cdot}(e_j)} (\nn X_{t})^{-1}_{ji}\Big\}\bigg],\end{split}\end{equation}
 where
 $$\dd (h_{t,\cdot}(e_j)):= \int_0^t \big\<\pp_s h_{t,s}(e_j),\d W_s\big\>$$ is the Malliavin divergence of $h_{t,\cdot}(e_j)$.
 Consequently
 \beq\label{Q6}
 \big|\E(\nn_{e_i} f)(X_t)\big]\big|\le \big(\E |f(X_t)|^p\big)^{\ff 1 p} \big(\E[|M_t(e_i)|^{\ff p{p-1}}\big]\big)^{\ff{p-1}p}
 \end{equation}
for $t\in (0,T] $ and $1\le i\le d_1+d_2 $.

By  \eqref{Q2} and \eqref{Q4}, there is a constant $c_4>0$ such that
\beq\label{Q7}
\big(\E[|M_t(e_i)|^{\ff p{p-1}}\big]\big)^{\ff{p-1}p}\le c_4 \sum_{j=1}^{d_1+d_2} 1_{\{\|(\nn X_t)^{-1}_{ji}\|_\infty>0\}}\bigg\{\E\bigg(\int_0^t |\pp_s h_{t,s}(e_j)|^2 \d s \bigg)^{\ff p{2(p-1)}} \bigg\}^{\ff{p-1}p}
\end{equation}
for any $ t\in (0,T]$ and $ 1\le i\le d_1+d_2$.

By \eqref{Q}, we have $\|Q_{t,s}^{-1}\|\le c_0^{-1} t s^{-2(k+1)}.$ Combining this with \eqref{Q0}, we may find a constant $c_5>0$ such that
\beg{align*}
&|\aa_{t,s}(e_j)|\le c_5 t^{-2k} + c_5 1_{\{j\le d_1\}} t^{-2k-1},\\
&|\pp_s \aa_{t,s}(e_j)|\le c_5 t^{-2k-1} + c_5 1_{\{j\le d_1\}} t^{-2k-2},\\
&|g_{t,s}(e_j)|\le c_5 t  + c_5 1_{\{j\le d_1\}}\ \ \text{ for } 0\le s<t\le T \text{ and } 1\le j\le d_1+d_2.
\end{align*}
Now noting that $\|\si(s)^{-1}\|\le K$, together with the previous estimates, we may conclude that there is  a constant $c_6>0$ such that
\beg{align*}
\pp_s h_{t,s}(e_j)|&=\big|\si(s)^{-1} \{\nn_{g_{t,s}(e_j), \aa_{t,s}(e_j)} b(s,X_s) -\pp_s \aa_{t,s}(e_j)\}\big|\\
& \le c_6 t^{-2k-1}+ c_6 1_{\{j\le d_1\}} t^{-2k-2}
\end{align*}
for any $ 0\le s<t\le T$ and for $1\le j\le d_1+d_2$.
This together with \eqref{Q7} enables us to find a constant $c_7>0$ such that
$$  \big(\E[|M_t(e_i)|^{\ff p{p-1}}\big]\big)^{\ff{p-1}p}\le c_7 \beg{cases}   t^{-2k-\ff 3 2},\ &\text{if}\ \sup_{j\le d_1} \|(\nn X_t)^{-1}_{ji}\|_\infty >0,\\
t^{-2k-\ff 1 2},\ &\text{otherwise}.\end{cases}$$
Combining this with \eqref{Q6} we derive \eqref{ES0} for some constant $c(p)>0$.

(d) For the case where $Z (s,x)= (s,x^{(2)})$ is independent of $x^{(1)}$, we have $\nn_j X_t^i=0$ for $i\ge d_1+1$ and $j\le d_1$, so that the previous estimate implies that
$$ \big(\E[|M_t(e_i)|^{\ff p{p-1}}\big]\big)^{\ff{p-1}p}\le c_7 t^{-2k-\ff 1 2}\ \ \ \forall  t\in (0,T], $$
where $ d_1+1\le i\le d_1+d_2$.
Combining this with \eqref{Q6} we derive we derive \eqref{ES} with some constant $c(p)>0$ and $\xi_t:= t^{-2k-\ff 12}$.
\end{proof}

\beg{lem}\label{LNN}  Let    $\eqref{R} $ and $ \eqref{R'}$ hold,  let $b\in C_b^1$, and  let  $Z$ be locally bounded  having  bounded $\nn Z$.
Then for any $p>1$ there exists a constant $c(p)>0$ such that the semigroup $P_t$ associated with \eqref{DS} satisfies the Harnack inequality
\beq\label{HII} |P_t f(x)|^p(x)\le \big(P_t|f|^p(y)\big)\e^{\ff{c(p)|x-y|^2}{t^{4k+3}}},\ \ \ t\in (0,T], x,y\in \R^{d_1+d_2}, f\in \B_b(\R^{d_1+d_2}).\end{equation}
\end{lem}

\beg{proof} (a) Let $\bar P_t$ be the Markov semigroup associated with \eqref{DS} for $b=0$.  By \cite[Corollary 4.3(1)]{WZ13} for $l_1=0$,   we find a constant $c_1(p)>0$ such that
\beq\label{A1} \hat P_t |f|(x)\le (\hat P_t |f|^{p^{\ff 1 3}} (y))^{p^{-\ff 1 3}}\e^{\ff{c_1(p)|x-y|^2}{t^{4k+3}}},\ \ \ t\in (0,T], x,y\in\R^{d_1+d_2}\end{equation}
holds for all $f\in \B_b(\R^{d_1+d_2}).$

On the other hand,  since $b$ is bounded,   there exists a constant $c_2(p)>0$ such that
$$P_t |f|\le \e^{c_2(p) t} (\hat P_t|f|^{p^{\ff 1 3}})^{p^{-\ff 1 3}}),\ \ \  \hat P_t |f|\le \e^{c_2(p) t} (P_t|f|^{p^{\ff 1 3}})^{p^{-\ff 1 3}}),\ \ t\in [0,T].$$
Combining this with \eqref{A1} we find a constant $c_3(p)>0$ such that
\beq\label{A2}   P_t |f|(x)\le (  P_t |f|^{p } (y))^{\ff 1 p}  \e^{c_3(p)+ \ff{ c_3(p)|x-y|^2}{t^{4k+3}}},\ \ \ t\in (0,T], x,y\in\R^{d_1+d_2}\end{equation}
holds for all $f\in \B_b(\R^{d_1+d_2}).$

Finally, since $\nn b$ and $\nn Z$ are bounded, $(\nn X_t)_{t\in [0,T]}$ is bounded as well. So,  there exists a constant $c_4>0$ such that
$$|\nn P_t f|\le c_4 P_t |\nn f|,\ \ \ t\in [0,T], f\in C_b^1(\R^{d_1+d_2}).$$
According to the proof of \cite[Theorem 2.2]{Ren}, this together  with \eqref{A2} implies \eqref{HII} for some constant $c(p)>0.$
\end{proof}

\section{Distribution dependent stochastic Hamilton system}

Consider  the following distribution dependent SDEs
\beq\label{DIS}
\beg{cases}
\d X_t^{(1)}= \big\{A X_t^{(1)} + B X_t^{(2)}+b(X_t, \L_{X_t})\big\} \d t,\\
\d X_t^{(2)}=  Z (t,X_t,\L_{X_t}) \d t
+ \si(t,\L_{X_t})\d W_t
\end{cases}
\end{equation}
 for $t\in [0,T]$, where  $X_t =(X_t^{(1)}, X_t^{(2)}) $   is $\R^{d_1}\times \R^{d_2}$ valued process. The coefficients  $A$, $B$,  $b, Z $ and $\si$ satisfy  the following assumption.

\beg{enumerate} \item[$(C_1)$]
$A, B$ and $b$ satisfy  conditions $1)$ and $2)$ in $(B_3),$      $Z (t,x,\mu)$ is differentiable in $x\in \R^{d_1+d_2}$, and  there exists a constant $K>0$ such that
$$
\|\nn b(t,\cdot,\mu)(x)-\nn b(t,\cdot,\mu)(y)\|\le K|x-y|,
$$
$$
|b(t,x,\mu)-b(t,y,\nu)|+  \|\si(t,\mu)-\si(t,\nu)\|\le K\big\{|x-y|+\W_2(\mu,\nu)\big\}
$$
$$
\|Z (t,0,\dd_0)|+ \|\si(t,\mu)\|+\|\si(t,\mu)^{-1}\|\le K
$$
for $t\in [0,T] $,  $x,y\in\R^{d_1+d_2} $ and $\mu,\nu\in \scr P_2(\R^{d_1+d_2})$.
%\beg{align*}& \\
%&\\
%&\end{align*}
\end{enumerate}

By, for instance, \cite[Theorem 2.1]{W18}, under this assumption  the SDE \eqref{DIS} is well-posed for distributions in $\scr P_2(\R^{d_1+d_2}),$ and $P_t^*\mu:= \L_{X_t}$ for the solution $X_t$ with initial distribution $\mu$ satisfies
\beq\label{WS}  \sup_{t\in [0,T]} \W_2(P_t^*\mu, P_t^*\nu) \le C \W_2(\mu,\nu),\ \ \forall \mu,\nu\in \scr P_2(\R^{d_1+d_2})\end{equation}
for some constant $C>0$.

\beg{thm}\label{T5} Assume that condition  $(C_1)$ is satisfied.
\beg{enumerate}
\item[$(1)$]
There exists a constant $c>0$ such that
\beq\label{EN1}
\Ent(P_t^*\mu|P_t^*\nu)\le  \ff c {t^{(4k+2)(4k+3)}}  \W_2(\mu,\nu)^2,\ \ \ \ \forall t\in (0,T].
\end{equation}
If $Z(t,x,\mu)=Z(t,x^{(2)},\mu)$ does not dependent on $x^{(1)}$, then
\beq\label{EN1*}
\Ent(P_t^*\mu|P_t^*\nu)\le   \ff c {t^{(4k+1)(4k+3)}}  \W_2(\mu,\nu)^2,\ \ \ \ \forall t\in (0,T].
\end{equation}

 \item[$(2)$]
If $Z (t,x,\mu)= Z(x,\mu)$ and $\si(t,\mu)=\si(\mu)$ do not depend on $t$, and there exist constants $c',\ll>0$ such that
$$
\W_2(P_t^*\mu,P_t^*\nu)^2\le c'\e^{-\ll t} \W_2(\mu,\nu)^2,
\  \ \forall t\ge 0 \ \mbox{ and }\ \ \forall  \mu,\nu\in \scr P_2(\R^{d_1+d_2}),
$$
then $P_t^*$ has a unique invariant probability measure $\bar\mu\in \scr P_2(\R^{d_1+d_2})$, and
$$
\Ent(P_t^*\mu|\bar\mu)\le   cc'\e^{-\ll(t-1)} \W_2(\mu,\bar\mu)^2
$$
for any $t\ge 0$ and for every $\mu\in \scr P_2(\R^{d_1+d_2})$.
\end{enumerate}
\end{thm}

\beg{proof} It suffices to prove the first assertion. To this end, given $(\mu,\nu\in \scr P_2(\R^{d_1+d_2}),$ let
\beg{align*}&Z_1^{(2)}(t,x):= Z(t,x,P_t^*\mu),\ \ Z_2^{(2)} (t,x):= Z(t,x,P_t^*\nu),\\
&  \si_1(t):= \si(t,P_t^*\mu)\ \ \  \si_2(t):=\si(t,P_t^*\nu),\ \ \  t\in [0,T].\end{align*}
Then  the desired estimates in Theorem \ref{T5}(1) follow from Corollary \ref{C3.1} and \eqref{WS}.
\end{proof}

To illustrate this result, we consider the following  typical example for $d_1=d_2=d$:
\beq\label{E21} \beg{cases} \d X_t^{(1)}= \big\{BX_t^{(2)} +b(X_t)\big\}\d t,\\
\d X_t^{(2)}=   \si(\L_{X_t}) \d W_t -
\left( B^* \nn V(\cdot,\L_{X_t})(X_t) + \bb B^* (BB^*)^{-1}X_t^{(1)}+X_t^{(2)}\right)
\d t,\end{cases}\end{equation}
where $\bb>0$ is a constant, $B$ is an  invertible $d\times d$-matrix, and
$$V: \R^{d}\times \scr P_2(\R^{2d})\to \R^{d}$$ is measurable and differentiable in $x^{(1)}\in\R^{d}$.
 Let
\beg{align*}
& {\psi}(x,y) :=\ss{|x^{(1)}-y^{(1)}|^2 +|B(x^{(2)}-y^{(2)})|^2}\ \ \text{ for } x,y\in \R^{2d},\\
&\W_2^{\psi}(\mu,\nu):=\inf_{\pi\in \C(\mu,\nu)} \bigg( \int_{\R^{2d}\times\R^{2d}} {\psi}^2 \d\pi\bigg)^{\ff 1 2}\ \ \ \text{ for }\mu,\nu\in\scr P_2(\R^{2d}).
\end{align*}

We assume that the following technical condition is satisfied.

\beg{enumerate} \item[$(C_2)$]   $V(\cdot,\mu)$ is differentiable such that $\nn V(\cdot,\mu)(x^{(1)})$ is Lipschitz continuous in $(x^{(1)},\mu)\in \R^{d}\times \scr P_2(\R^{2d}).$
Moreover, there exist   constants $\theta_1, \theta_2\in \R$ with
$$
\theta_1+\theta_2<\bb,
$$  such that
 \beg{align*}
 & \big\< BB^* \{\nn V(\cdot,\mu)(x^{(1)})-\nn V(\cdot,\nu)(y^{(1)})\}, \  x^{(1)}-y^{(1)}+(1+\bb)B(x^{(2)}-y^{(2)})\big\>\\
&-\ff {1+\bb}{2\bb} \|B\{\si(\mu)-\si(\nu)\}\|_{HS}^2\ge -\theta_1{\psi} (x,y)^2 -\theta_2 \W_2^{\psi} (\mu,\nu)^2
\end{align*}
for any $x,y\in \R^{2d}$ and $\mu,\nu\in \scr P_2(\R^{2d})$.
\end{enumerate}

\beg{cor}\label{CC}
Assume that condition $(C_2)$ is satisfied. Let
\beq\label{KK} \kk:=\ff{ 2(\bb-\theta_1-\theta_2)}{2+2\bb+\bb^2+\ss{\bb^4+4}}>0.\end{equation}
 For any $\kk'\in (0,\kk)$, when $\|\nn b\|_\infty$ is small enough,   $P_t^*$ has a unique invariant probability measure $\bar \mu \in \scr P_2(\R^{2d}),$  and there exists a constant   $c >0$ such that
\beq\label{KK0}
\W_2(P_t^*\mu,\bar\mu)^2 + \Ent(P_t^*\mu|\bar \mu)  \le  \ff{c\e^{-2\kk' t}}{(1\land t)^{3}} \W_2(\mu,\bar\mu)^2
 \end{equation}
 for any $t>0$ and $ \mu\in \scr P_2(\R^{2d})$.

\end{cor}

\beg{proof} The proof is completely similar to that of \cite[Lemma 5.2]{RW21} where $\si(\mu)=\si$ does not depend on $\mu$. By Theorem \ref{T5}, it suffices to find a constant $c'>9$ such that
\beq\label{EII}
\W_2(P_t^*\mu,P_t^*\nu)^2\le c'\e^{-2\kk t} \W_2(\mu, \nu)^2
\end{equation}
for any $ t>0$ and $\mu,\nu\in \scr P_2(\R^{2d}) $.

a) Let
 \beq\label{AC1} a:=\Big(\ff{1+\bb+\bb^2}{1+\bb}\Big)^{\ff 1 2} ,\ \ r:=  a -\ff{\bb}{ a} =\ff 1 {\ss{(1+\bb)(1+\bb+\bb^2)}}\in (0,1).\end{equation}
 Define  the distance
\beq\label{AC2}  \bar \psi(x,y):=\ss{a^2|x^{(1)}-y^{(1)}|^2 +|B(x^{(2)} -y^{(2)})|^2+ 2 r a \<x^{(1)}-y^{(1)}, B(x^{(2)}-y^{(2)})\>}.
 \end{equation}   According to the proof of  \cite[Lemma 5.2]{RW21}, we have
 \beq\label{AC3}   \bar\psi (x,y)^2\le \ff{2+2\bb+\bb^2 +\ss{\bb^4+4}}{2(1+\bb)}   \psi (x,y)^2, \ \ \forall x,y\in\R^{2d},\end{equation}
 and there exists a constant $C>1$ such that
 \beq\label{ACC}
 C^{-1}|x-y|\le  \bar\psi (x,y) \le C|x-y|,
 \ \ \forall x,y\in\R^{2d}.
 \end{equation}

  b) Let  $X_t$ and $Y_t$ solve \eqref{E21} with $\L_{X_0}=\mu, \L_{ Y_0}=\nu$ such that
\beq\label{AC4} \W_2(\mu,\nu)^2 =\E\big[|X_0- Y_0|^2\big].\end{equation}
Let $\Xi_t=X_t-Y_t$, $\mu_t=P_t^*\mu:=\L_{X_t}$ and $  \nu_t:=P_t^*\nu=\L_{Y_t}$.  By using $(C_2)$,     It\^o's formula, and noting that \eqref{AC1} implies
$$
a^2-\bb-ra=0,\ \ 1-ra=ra\bb =\ff\bb {1+\bb},
$$
we obtain
\beg{align*} &\ff 1 2 \d \left(\bar\psi(X_t,Y_t)^2\right)=\ff 1 2 \|B\left(\si(\mu_t)-\si(\nu_t)\right)\|_{HS}^2
 + \big\<a^2\Xi_t^{(1)} + r a B\Xi_t^{(2)}, B\Xi_t^{(2)} + b(X_t)-b(Y_t)\big\>\d t  \\
&\quad  -  \big\<B^*B \Xi_t^{(2)} + r a B^*\Xi_t^{(1)},\ \bb B^*(BB^*)^{-1} \Xi_t^{(1)}  +  \Xi_t^{(2)}  \big\>\d t \\
&\quad+\big\<B^*B \Xi_t^{(2)} + r a B^*\Xi_t^{(1)},\ B^*\{\nn^{(1)} V(Y_t^{(1)},  \nu_t) -\nn^{(1)} V(X_t^{(1)},\mu_t)\}\big\> \d t\\
&\le \Big\{-(1-ra)|B\Xi_t^{(2)} +(a^2-\bb-ra) \<\Xi_t^{(1)}, B\Xi_t^{(2)} \> +\big[\|\nn b\|_\infty (a^2+r a)- ra\bb \big]     |\Xi_t^{(1)}|^2 \\
&\quad + \big\<B^*B \Xi_t^{(2)} + (1+\bb)^{-1}B^*\Xi_t^{(1)},\ B^*\{\nn^{(1)} V(Y_t^{(1)},  \nu_t) -\nn^{(1)} V(X_t^{(1)},\mu_t)\} \big\>\Big\}\d t\\
&\le \Big\{\ff{\theta_2}{1+\bb}\W_2^{\psi}(\mu_t,\nu_t)^2-\ff{\bb-\theta_1}{1+\bb}{\psi}(X_t,Y_t)^2 + \|\nn b\|_\infty (a^2+r a)     |\Xi_t^{(1)}|^2\Big\}\d t.\end{align*}
By \eqref{AC3} and the fact   that
$$\W_2^{\psi}(\mu_t,\nu_t)^2\le \E[{\psi}(X_t,Y_t)^2],$$
for $\kk>0$ in \eqref{KK},  when $\|\nn b\|_\infty$ is small enough we find a constant $\kk'\in (0,\kk)$ such that
we obtain
\beg{align*}
&\ff 1 2 \left(\E[ \bar\psi(X_t,Y_t)^2]- \E[ \bar\psi(X_s,Y_s)^2]\right)\\
&\le \|\nn b\|_\infty (a^2+r a)    \int_s^t  \E[|\Xi_u^{(1)}|^2] \d u -\ff{\bb-\theta_1-\theta_2}{1+\bb} \int_s^t \E [{\psi}(X_u,Y_u)^2]\d u \\
&\le -\kk'  \int_s^t \E [ \bar\psi(X_u,Y_u)^2]\d u,\ \ t\ge s\ge 0.
\end{align*}
By  Gronwall's inequality, we then deduce that
$$
\E[ \bar\psi(X_t,Y_t)^2]\le \e^{-2\kk' t} \E[ \bar\psi(X_0,Y_0)^2]
$$
for $t\ge 0$.
Combining this with \eqref{ACC} and \eqref{AC4}, we may conclude that there is  a constant  $c>0$ such that \eqref{EII} holds.
\end{proof}

To conclude this paper, we present the following example of degenerate nonlinear granular media equations, see \cite{CMV} and \cite{GLW} for the study of non-degenerate linear granular media equations.

\paragraph{Example 4.1 (Degenerate nonlinear  granular media equation).}  Let $ d\in \mathbb N$ and $W\in C^\infty(\R^d\times\R^{2d}).$ Consider  the following PDE for probability density functions $(\rr_t)_{t\ge 0}$ on $\R^{2d}=\R^d\times\R^d$:
\beq\label{*EN2} \beg{split} &\pp_t\rr_t(x)= \ff 1 2 {\rm tr}\big\{\si(\rr_t)\si(\rr_t)^*(\nn^{(2)})^2\big\} \rr_t(x) -\<\nn^{(1)}\rr_t(x), x^{(2)}+b(x)\>\\
&+ \<\nn^{(2)} \rr_t(x), \nn^{(1)}(W\circledast\rr_t)(x^{(1)}) + \bb x^{(1)}+x^{(2)}\>,
\end{split}
\end{equation}
where $x=(x^{(1)},x^{(2)})\in\R^{2d} $, $t\ge 0 $.  $\bb>0$ is a constant,     and
$$(W\circledast \rr_t)(x^{(1)}):=\int_{\R^{2m}} W(x^{(1)},z)\rr_t(z) \d z,\ \ x^{(1)}\in \R^d$$  stands for the mean field interaction.

  If there exist constants $\theta,\aa >0$ with
  $$
  \theta\Big(\ff 1 2+\ss{2+2\bb+\bb^2}\Big)+\ff{\aa(1+\bb)}{2\bb} <\bb,
  $$
  such that
\beq\label{*EN3}
\beg{split}
&|\nn W(\cdot,z)(v)-\nn W(\cdot,\bar z)(\bar v) |\le \theta \big(|v-\bar v|+|z-\bar z|\big),\ \ \forall v,\bar v\in \R^d, \text{ and } \forall z,\bar z\in \R^{2d},\\
&\|\si(\mu)-\si(\nu)\|_{HS}^2 \le \aa\W_2(\mu,\nu)^2,\ \ \ \forall \mu,\nu\in\scr P_2(\R^{2d}),
\end{split}
\end{equation}
then  for any $\kk'\in (0,\kk)$, when  $\|\nn b\|_\infty$ is small enough   there exists a unique probability measure  $\bar \mu \in\scr P_2( \R^{2d})$   and a constant $c>0$ such that for any probability density functions $(\rr_t)_{t\ge 0}$ solving \eqref{*EN2},  $\mu_t(\d x):=\rr_t(x)\d x$ satisfies \beq\label{ERR}
\W_2(\mu_t,\bar \mu)^2+\Ent(\mu_t|\bar\mu)     \le c\e^{-\kk' t}  \W_2(\mu_0,\bar \mu)^2,\ \ \forall t\ge 1
\end{equation}
where
$$\kk= \ff{ 2\bb-\theta-2\theta \ss{2+2\bb+\bb^2}-\aa(1+\bb^{-1})}{2+2\bb+\bb^2+\ss{\bb^4+4}} >0.$$

To prove this claim, let   $(X_t,Y_t)$ solve \eqref{E21} for
\beq\label{GPP} B:= I_{d},\ \ \psi(x,y)=|x-y|,\ \text{ and } V(x, \mu):=  \int_{\R^{2d}} W(x,z)\mu(\d z).\end{equation}
As shown in the proof of \cite[Example 2.2]{RW21},   $\rr_t$ solves \eqref{*EN2} if and only if $\rr_t(x)=\ff{\d (P_t^*\mu)(\d x)}{\d x}$,  where $P_t^*\mu:=\L_{X_t}.$

 By Corollary \ref{CC},  we only need to verify  $(C_2)$  for $B, V$ in \eqref{GPP}  and
\beq\label{TTH} \theta_1= \theta\Big(\ff 1 2 +\ss{2+2\bb+\bb^2}\Big),\ \ \theta_2= \ff \theta 2 \ss{2+2\bb+\bb^2}+\ff{\aa(\bb+1)}{2\bb},\end{equation}
so that the desired assertion holds for
$$ \kk:=\ff{ 2(\bb-\theta_1-\theta_2)}{2+2\bb+\bb^2+\ss{\bb^4+4}}>0.$$

For simplicity, let $\nn^v$ denote the gradient in $v$.
By \eqref{*EN3} and $V(x,\mu):=\mu(W(x,\cdot))$, for any constants $\aa_1,\aa_2,\aa_3>0$ we have
\beg{align*} I &:= \big\<\nn^{x^{(1)}} V(x^{(1)},\mu)- \nn^{y^{(1)}} V(y^{(1)},\nu), x^{(1)}- y^{(1)} +(1+\bb)(x^{(2)}-  y^{(2)})\big\> \\
&\le   \int_{\R^{2m}} \big\< \nn^{x^{(1)}} W(x^{(1)}, z) - \nn^{y^{(1)}} W(y^{(1)}, z), \ x^{(1)}-y^{(1)} +(1+\bb)(x^{(2)}-y^{(2)})\big\>\mu(\d z) \\
&\qquad +\big\<\mu(\nn^{y^{(1)}}  W(y^{(1)},\cdot))- \nu(\nn_{y^{(1)}} W(y^{(1)},\cdot)),  x^{(1)}-y^{(1)} +(1+\bb)(x^{(2)}-y^{(2)})\big\>\\
&\ge  -\theta\big\{ |x^{(1)}-y^{(1)}|+ \W_1(\mu,\nu)\big\} \cdot\big(|x^{(1)}-y^{(1)}|+(1+\bb)|x^{(2)}-y^{(2)}|\big)\\
&\ge - \theta (\aa_2+\aa_3) \W_2(\mu,\nu)^2\\
&\quad- \theta\Big\{\Big(1+\aa_1+\ff 1 {4\aa_2}\Big)|x^{(1)}-y^{(1)}|^2 +(1+\bb)^2 \Big(\ff 1 {4 \aa_1}+\ff 1 {4\aa_3}\Big)|x^{(2)}-y^{(2)}|^2\Big\}.\end{align*}
 Take
$$\aa_1= \ff{\ss{2+2\bb+\bb^2}-1} 2,\ \ \aa_2= \ff 1 {2\ss{2+2\bb+\bb^2 }},\ \ \text{ and } \aa_3=\ff {(1+\bb)^2}{2\ss{2+2\bb+\bb^2 }}.$$
We have
\beg{align*}&1+\aa_1+\ff 1 {4\aa_2} =\ff 1 2+ \ss{2+2\bb+\bb^2},\\
&(1+\bb)^2\Big(\ff 1 {4 \aa_1}+\ff 1 {4\aa_3}\Big)=\ff 1 2+ \ss{2+2\bb+\bb^2},\\
& \ \aa_2+\aa_3= \ff 1 2 \ss{2+2\bb+\bb^2}.\end{align*}
Combining this with \eqref{*EN3} and \eqref{TTH}, we derive
 $$I-\ff{\bb+1}{2\bb}\|\si(\mu)-\si(\nu)\|_{HS}^2\ge  -\theta_1 |x-y)|^2-\theta_2 \W_2(\mu,\nu)^2,$$
and therefore condition  $(C_2)$  is satisfied for $B,\psi$ and $V$ in \eqref{GPP} .

\end{document}